\documentclass[12pt,twoside, a4paper]{article}

\usepackage{latexsym} 
\usepackage{a4}

\newtheorem{bsp}{Example} 
\newtheorem{Def}{Definition} 
\newtheorem{satz}{Proposition} 
\newtheorem{folg}[satz]{Corollary}
\newtheorem{theo}[satz]{Theorem}
\newtheorem{lem}[satz]{Lemma}

\font\fr=eufm10

\def\bbm#1{\mbox{\boldmath$#1$}}

\def\ffr#1{\mbox{\fr #1}}

\def\pro{P\,r\,o\,o\,f.\enspace}

\def\bst{$\Box$}

\let\empt=\emptyset

\def\fb{\unskip\nobreak\vadjust{\kern 2pt}\hfil\break
 \parfillskip=0pt\vadjust{}\hfill{}}

\def\Eend{\mathop{\rm End}\nolimits}

\def\dim{\mathop{\rm dim}\nolimits}

\def\id{\mathop{\rm id}\nolimits}

\def\ii{\mathop{\rm i}\nolimits}

\def\begeq{\begin{equation}}
\def\eneq{\end{equation}}
\def\begbsp{\begin{bsp}\hspace{-1.2ex}. \rm}
\def\enbsp{\bst\end{bsp}}
\def\begdef{\begin{Def}\hspace{-1.2ex}.\hspace{1ex} \rm} 
\def\endef{\bst\end{Def}}
\def\begsatz{\begin{satz}\hspace{-1.2ex}.\hspace{1ex}}
\def\ensatz{\end{satz}}
\def\begfolg{\begin{folg}\hspace{-1.2ex}.\hspace{1ex}}
\def\enfolg{\end{folg}}
\def\begub{\small\begin{ub}\hspace{-1.2ex}.\ \rm}
\def\enub{\end{ub}\normalsize}
\def\beglem{\begin{lem}\hspace{-1.2ex}.\hspace{1ex}}
\def\enlem{\end{lem}}

\begin{document}

\title{M\"obius invariants for pairs of spheres $(S_1^m, S_2^l)$ in the
M\"obius space $S^n$}

\author{Rolf Sulanke}


\date{February 22, 1999}

\maketitle

\begin{abstract}
In this article we construct a complete system of
M\"obius-geometric invariants for pairs $(S^m, S^l),\,l \leq  m,$
of spheres contained in the M\"obius space $S^n$. It consists of
$n-m$ generalised stationary angles. We interpret these
invariants geometrically. \footnote{Preprint. Supported by the Institute of
Mathematics,Humboldt University, Berlin.}
\end{abstract}
\section{Introduction}
The $n$-dimensional real \em M\"obius geometry \em is the
$n$-dimensional sphere $S^n$ considered under the action of the
\em M\"obius group \em $G_n$, which is the group of all conformal
transformations of the Riemannian sphere $S^n$ with the standard
metric of constant sectional curvature. By  {\sc Liouville}'s
theorem  the group $G_n$ coincides with the group of those
diffeomorphisms of $S^n$ which transform hyperspheres of the \em
M\"obius space \em $S^n$ into hyperspheres, or generally
$m$-spheres into $m$-spheres, see e. g. {\sc W.~Blaschke,
G.~Thomsen} \cite{BT:29}, or {\sc M.~A.~Akivis, V.~V.~Goldberg}
\cite{AG:96}.  If one considers $S^n$ as a hyperellipsoid in the
$n+1$-dimensional real projective space $\bbm P^{n+1}$, the
$m$-spheres of $S^n$ are the non empty intersections of $S^n$
with projective, not tangent $m+1$-planes. Therefore, the group
$G_n$ consists of all projective transformations preserving the
hyperellipsoid $S^n$. The M\"obius group $G_n$ appears as the
group of  projective transformations generated by
pseudo-orthogonal transformations of the underlying vector space
provided with a scalar product of index 1 (the $n+2$-dimensional
{\sc Minkowski}-space).
In  {\sc B. A. Rosenfeld}'s book \cite{Ro:x} one finds the
M\"obius-geometric  definition of stationary angles for pairs of
subspheres having equal dimension $m = l$. We refine his method
to handle the general case, and carry out the classification. In
spite of the fact that real M\"obius geometry is based on
pseudo-euclidean linear algebra we obtain the classification only
using  the elementary theory of selfadjoint operators in
finite-dimensional euclidean vector spaces. Another approach to
the same problem is contained in the paper \cite{AFC:98} of {\sc
A.~Montesinos Amilibia, M.~C.~Romero Fuster} and {\sc E.~Sanabria
Codesal}. I thank these colleagues very much for the preprint,
for interesting discussions, and their kind support of my visit
in Valencia.  

In section 2 we describe the mentioned above projective model of
the M\"obius geometry and fix the basic notations. Section 3
contains the case of hyperspheres, for which a generalised angle
is the only conformal invariant. In Section 4 the (generalised)
stationary angles between subspheres of arbitrary, in general
distinct dimensions are introduced, and  the classification of
pairs of such subspheres is carried out.

\section{The projective model of the M\"obius space}

Let $\bbm{V}$ denote the \em $(n+2)$-dimensional real vector
space\em, and $<\ffr x, \ffr y>,\, \ffr x, \ffr y \in \bbm{V},$ a
given non-degenerate symmetric bilinear form of index one. We
call it \em the scalar product\em. A basis $(\ffr b_i),\, i =
1,\ldots, n+2,$ is said to be \em orthogonal\em,\, if the scalar
products fulfil
\begeq
<\ffr b_i, \ffr b_j> = \delta_{ij}\epsilon_i,\, i, j = 1,\ldots,
n+2,
\eneq 
\begeq
\epsilon_i = 1 \mbox{ for } i= 1,\ldots, n+1,\, \epsilon_{n+2} =
- 1.
\eneq 
In (1) $\delta_{ij}$ denotes the Kronecker-Symbol. The \em
projective space \em belonging to $\bbm{V}$ is denoted by
$\bbm{P}$,  as usual it is defined to be the set of all
one-dimensional subspaces of $\bbm{V}$. For $\ffr x\in \bbm{V},
\ffr x\neq \ffr o,$ we denote by
\begeq
x := [\ffr x] := \mbox{ linear hull}\{\ffr x\} \in \bbm{P} 
\eneq 
the corresponding point of the projective space $\bbm{P}$. Any
subspace $\bbm W$ of dimension $\dim \bbm W = m+1$ is at the same
time a \em projective $m$-plane\em, which is defined as the set
of points obtained by (3) for $\ffr x\in W$. The vector
coordinates of $\ffr x \in \bbm{V}$ with respect to any basis of
$\bbm{V}$ are the corresponding \em homogeneous coordinates\em\
of the point $x=[\ffr x] \in \bbm{P}$; they are defined by $x$ up
to a common non-zero factor. The \em $n$-dimensional sphere\em\
$S^n$ is defined to be the projective hyperquadric
\begeq
x\in S^n : \Longleftrightarrow x=[\ffr x]\in \bbm{P} \mbox{ and }
<\ffr x, \ffr x> = 0.
\eneq 
With respect to any orthogonal basis we get from (4) the equation
of $S^n$ in homogeneous coordinates $x_i$ of $x$ in the normal
form
\begeq
\sum_{i=1}^{n+1} x_i^2 - x_{n+2}^2 = 0.
\eneq 
Obviously, for any pseudo-orthogonal coordinate system and any
$x\in S^n$ one has $x_{n+2}\neq 0$. Therefore we may norm the
coordinates of the points of $S^n$ by $\xi_i := x_i/x_{n+2}$, and
obtain the usual equation of the unit hypersphere:
\begeq
\sum_{i=1}^{n+1} \xi_i^2 = 1.
\eneq 
The isotropy group $G_n$ of $S^n$, that is the subgroup of all
projective transformations of $\bbm{P}$ preserving $S^n$, is
called the \em M\"obius group\em. It is well known, and easy to
prove,  that $G_n$ is isomorphic to the pseudo-orthogonal group
acting on $\bbm V$, factored by $\{\bf 1, -\bf 1\}$,
where $\bf 1$ denotes the identity transformation of $\bbm V$:
\begeq
G_n \cong {\bf O}(n+1,1)/\{{\bf 1}, -{\bf 1}\}\cong {\bf
O}(n+1,1)_0\cup{\bf O}(n+1,1)_1.
\eneq 
Here ${\bf O}(n+1,1)_0$ is the component of the identity of
${\bf O}(n+1,1)$, and ${\bf O}(n+1,1)_1$ is the coset of the
transformations reversing only the orientation of space (or of
time, but not of both) in the relativistic interpretation of
$[\bbm{V}, <\, ,\, >]$. Therefore also $G_n$  has two components,
where the component $SG_n$ of the identity preserves the
orientation, and the other transformations reverse the
orientation of $S^n$. The sphere $S^n$ considered under the
projective action of the group $G_n$ described above, is called
the \em $n$-dimensional M\"obius space\em, and the geometry of
this transformation group is the \em M\"obius geometry\em. 

In spite of the fact that the M\"obius geometry can be developed
in an abstract manner considering only the geometry of objects
contained in $S^n$, we prefer the projective point of view, since
then all tools of the pseudo-euclidean linear algebra can be
applied directly. Clearly, the points $x\in S^n$ are represented
by \em isotropic vectors\em\ $\ffr x \neq \ffr o,\, <\ffr x, \ffr
x>=0, x = [\ffr x]$. With respect to their position to $S^n$,
there exist three types of projective subspaces:

A. The projective $m+1$-planes which \em intersect the M\"obius
space $S^n$ in an $m$-sphere; \em they correspond to the \em
pseudo-euclidean subspaces \em of dimension $m+2,\, m = 0,\ldots,
n$, and index $1$. A sphere $S^0$ is a set of two points, and for
$m = n$ we get the whole M\"obius space. 

B. The projective $m$-planes which \em are disjoint to the
M\"obius space $S^n$; \em they correspond to the \em euclidean
subspaces \em of dimension $m+1,\, m = 1,\ldots, n$, on which the
restriction of the scalar product is positive definite.

C. The projective $m$-planes which are \em tangent to the
M\"obius space $S^n$; \em they correspond to the \em isotropic
subspaces \em of dimension $m+1,\, m = 1, \ldots, n$, on which
the restriction of the scalar product degenerates and has rank
$m$. Any such $m$-plane contains exactly one point of the
M\"obius space $S^n$, the \em contact point.  \em 

The \em orthogonality \em $\bbm U\mapsto \bbm U^\bot$ with
respect to the pseudo-euclidean scalar product in $\bbm V$
defines an involutive bijection of the lattice of all subspaces
of $\bbm V$, whose projective interpretation is called a \em
polarity. \em It has the following properties:
\begin{eqnarray}
\dim \bbm U = k &\Longleftrightarrow &\dim \bbm U^\bot = n+2-k,\\
 \bbm U \subset \bbm W &\Longleftrightarrow &\bbm W^\bot \subset
\bbm U^\bot,\\
\bbm U \mbox{ pseudo-euclidean } &\Longleftrightarrow &\bbm
U^\bot \mbox{ euclidean},\\
\bbm U \mbox{ isotropic } &\Longleftrightarrow &\bbm U^\bot
\mbox{ isotropic},\\
\bbm U \mbox{ is not isotropic } &\Longleftrightarrow &\bbm U
\oplus \bbm U^\bot =\bbm V,\\
\bbm U \mbox{ is isotropic } &\Longleftrightarrow &\dim \bbm U
\cap \bbm U^\bot = 1,\\
(\bbm U + \bbm W)^\bot & = & \bbm U^\bot \cap \bbm W^\bot,\\
(\bbm U \cap \bbm W)^\bot & = &\bbm U^\bot + \bbm W^\bot.
\end{eqnarray} 

Of course, the case A of intersecting  projective $(m+1)$-planes
$\bbm U^{m+1}$ (projective dimension!)  leads to the natural
definition of an $m$-sphere $S^m_1 \subset S^n$ as the
intersection 
\begeq 
S^m_1 := \bbm U^{m+1}\cap S^n.
\eneq 
Since the scalar product restricted to the corresponding
$(m+2)$-dimensional vector subspace $\bbm U$ is pseudo-euclidean
of index 1 and rank $m+2$, the set $S^m_1$ is a hyperellipsoid in
$\bbm U$. Since the pseudo-orthogonal group acts transitively on
the set of all orthogonal frames one obtains as an immediate
consequence of properties (8) - (15):

\begsatz 
Definition \em (16) \em is an equivariant bijection between the
set of all $m$-spheres in $S^n$ and the set $X_m$ of all
$(m+2)$-dimensional pseudo-euclidean subspaces of the vector
space $\bbm V$, $ m = 0, 1,\ldots, n$. This set is an open
submanifold of the Grassmann manifold $G_{n+2, m+2}$ of all
$(m+2)$-dimensional subspaces of $\bbm V$; the group $O(n+1,1)$
acts transitively on $X_m$. Thus, the M\"obius group $G_n$ acts
transitively on the set of all $m$-spheres, which by \em (16) \em
can be identified with $X_m$. The orthogonality \em (10) \em
yields a bijection of $X_m$ with the set of all
$(n-m)$-dimensional euclidean subspaces of $\bbm V$, which is an
open submanifold of the Grassmann manifold of all
$(n-m)$-dimensional subspaces of $\bbm V$, and by \em (12) \em a
bijection of $X_m$ with all decompositions of $\bbm V$ into a
direct sum of an $(m+2)$-dimensional pseudo-euclidean subspace
and its orthogonal complement. \bst
\ensatz

In the following we identify $X_m$ with the set of all
$m$-spheres of $S^n$.

\section{Hyperspheres} 

By proposition 1 and (8) the set $X_{n-1}$ of all hyperspheres
$S^{n-1}\subset S^n$ can be identified  with the set of all
euclidean one-dimensional subspaces of $\bbm V$. Any such
subspace contains two spacelike unit vectors $\ffr n$; therefore
the hyper-hyperboloid $H^{n+1}$ defined by the equation
 \begeq
 <\ffr n, \ffr n> = n^2_1 + n^2_2+\ldots + n^2_{n+1} -n^2_{n+2} =
1
\eneq 
is a twofold covering of the manifold $X_{n-1}$ of non-oriented
hyperspheres. We denote the hypersphere corresponding to $\ffr n$
by 
\begeq
S(\ffr n) := S^n\cap [\ffr n]^\bot,\; \ffr n  \in H^{n+1},
\eneq 
and show
\begsatz 
Let $(N, M)$, $(\hat N,\, \hat M)$ be two pairs of non-oriented
distinct hyperspheres of $S^n,\, n\geq 1$:
$$ N = S(\ffr n),\, M = S(\ffr m),\, \hat N = S(\hat{\ffr n}),\,
\hat M = S(\hat{\ffr m}),\: \ffr n, \ffr m, \hat{\ffr n},
\hat{\ffr m}\in H^{n+1}.$$
Then $(N, M)$ is $G_n$-equivalent to $(\hat N,\, \hat M)$ if and
only if the absolute values of the scalar products coincide:
\begeq
|<\ffr n, \ffr m>| = |<\hat{\ffr n}, \hat{\ffr m}>|.
\eneq 
\ensatz
\pro For the proof we adapt an orthogonal frame $(\ffr b_i)$ to
the pair $(N,\, M)$, in which the position of $N, M$ depends on 
$|<\ffr n, \ffr m>|$ only; then the transformation defined by 
$$ g: g\ffr b_i = \hat{\ffr b}_i,\, i = 1, \ldots, n+2, $$
where $(\hat{\ffr b}_i)$ denotes the adapted frame corresponding
to $(\hat N, \hat M)$, is a M\"obius transformation with $(g N, g
M) = (\hat N,\, \hat M)$. To this aim we have to distinguish
three cases:

Case 1: $|<\ffr n, \ffr m>| < 1$, the subspace $[n, m]$ is \em
euclidean. \em Then we choose $\ffr b_1, \ffr b_2$ such that
$$ \ffr n =\ffr b_1,\: \ffr m = \ffr b_1 \cos \alpha + \ffr b_2
\sin \alpha, \mbox{ with } \cos \alpha = |<\ffr n, \ffr m>|,\, 0
< \alpha \leq \pi/2,$$
and complete $\ffr b_1, \ffr b_2$ to an orthogonal basis of $\bbm
V$. Doing the same with $(\hat N,\, \hat M)$, we obtain
$$ g \ffr n = \hat{\ffr n},\: g \ffr m = g(\ffr b_1 \cos \alpha +
\ffr b_2 \sin \alpha) = \hat{\ffr b}_1 \cos \alpha + \hat{\ffr
b}_2 \sin \alpha = \hat{\ffr m},$$
what was to show.

Case 2: $|<\ffr n, \ffr m>| > 1$, the subspace $[n, m]$ is \em
pseudo-euclidean. \em Then we choose $\ffr b_1, \ffr b_{n+2}$
such that
$$ \ffr n =\ffr b_1,\: \ffr m = \ffr b_1 \cosh \alpha + \ffr
b_{n+2} \sinh \alpha, \mbox{ with } \cosh \alpha = |<\ffr n, \ffr
m>|,\, 0 < \alpha,$$
and complete $\ffr b_1, \ffr b_{n+2}$ to an orthogonal basis of
$\bbm V$. 

Case 3: $|<\ffr n, \ffr m>| = 1$, the subspace $[n, m]$ is \em
isotropic. \em Then we choose $\ffr b_1, \ffr b_2, \ffr b_{n+2}$
such that
$$ \ffr n =\ffr b_1,\: \ffr m = \ffr b_1  + \ffr b_2 + \ffr
b_{n+2}, $$
and complete $\ffr b_1, \ffr b_2, \ffr b_{n+2}$ to an orthogonal
basis of $\bbm V$. Since $n\geq 1$, and the hyperspheres are
distinct, such a choice is always possible.

Concluding in cases 2, 3 as in case 1 one finishes the proof.
\bst

Now we interpret the three cases mentioned in the proof of
proposition 2 geometrically.

 Case 1. The hyperspheres $N, M$ intersect. The intersection is a
sphere of dimension $n-2$, and the cosinus of the intersection
angle is given by 
\begin{eqnarray}
\cos \alpha &=& |<\ffr n, \ffr m>|,\\ 
<\ffr n,\ffr m> &=& \frac{r^2 + R^2 - d^2}{2rR}, 
\end{eqnarray}
where $r, R$ are the radii of $N, M$, and $d$ is the distance of
the centres of these hyperspheres in the euclidean model of the
M\"obius space $S^n$, see below.

Case 2. The hyperspheres $N, M$ are disjoint. Formula (21)
remains valid; it expresses the conformal invariant of the two
hyperspheres by their euclidean invariants. Clearly, since each
conformal invariant is a fortiori an isometric invariant, it is
always possible to describe it in terms of a complete isometric
invariant system of the given objects. Using the equation $\cosh
\alpha = \cos \ii\alpha$ some authors introduce "imaginary
intersection angles", and interpret (20) in this sense, see e.g.
{\sc B.~A.~Rosenfeld} \cite{Ro:x}.

 Case 3. The hyperspheres $N, M$ have only one point in common,
therefore they are tangent to each other at this point.  

To prove these statements we denote by $\bbm W, \bbm U$ the
$n+1$-dimensional pseudo-euclidean subspaces of $\bbm V$ defining
the hyperspheres $N, M$. Since these hyperspheres are distinct,
it follows
$$ \dim(\bbm W + \bbm U) = n+2,\mbox{ and } \dim(\bbm W\cap \bbm
U) = n.$$
In case 1 the subspace spanned by $\ffr n, \ffr m$ 
$$[\ffr n, \ffr m] = \bbm W^\bot + \bbm U^\bot= (\bbm W \cap \bbm
U)^\bot$$
is euclidean, therefore the intersection $\bbm W \cap \bbm U$ is
pseudo-euclidean and defines an $n-2$-dimensional sphere. In case
2 this intersection is euclidean, and in case 3 isotropic, from
which the qualitative contents of the disjunction follows. Here
we applied formulas (10), (11), and (15).   It remains to prove
(20) and (21). If (21) is proved, then in case 1 formula (20) is
simply the cosinus theorem. Thus it suffices to show (21). To do
this we construct Riemannian and euclidean models of the M\"obius
space $S^n$ and derive a formula which relates the vectors $\ffr
n\in H^{n+1}$ with the centre and the radius of the hypersphere
$S(\ffr n)$ in the euclidean model. Let $(\ffr a_i), i =
1,\ldots, n+2,$ be a fixed orthogonal frame in the
pseudo-euclidean vector space $\bbm V$. The \em Riemannian model
\em is simply the unit n-sphere (6) with centre $\ffr o$ in the
euclidean subspace spanned by $\ffr a_1,\ldots, \ffr a_{n+1}$; it
can be considered also as the intersection of the isotropic cone
with the hyperplane $x_{n+2} = 1$, and therefore one has $\xi_i =
x_i, i = 1,\ldots, n+1$. The stereographical projection $sp$ of
$S^n$ from its north pole $\ffr a_{n+1}$ onto its equatorial
hyperplane $E^n$ spanned by the $\ffr a_i, i = 1,\ldots, \ffr
a_n$, is given by the formula
\begin{eqnarray}
x = \sum_{i= 1}^{n+1}\ffr a_i \xi_i\in S^n & \mapsto & sp(x) :=
\sum_{i= 1}^n\ffr a_i \xi_i/(1-\xi_{n+1})\in E^n, x\neq \ffr
a_{n+1},\\
sp(\ffr a_{n+1})&=& \infty.
\end{eqnarray} 
Of course, here (6) is supposed; the north pole $\xi_{n+1} = 1$
is
mapped onto the infinite point $\infty$, compactifying the $E^n$
to
$S^n$. Since the stereographical projection (22) is conformal,
and
transforms $k$-spheres into $k$-spheres, we consider the
compactified
$E^n$ as the \em euclidean model \em of the M\"obius geometry.
Hereby
one has to take into account that the $k$-planes, as the images
of the
$k$-spheres containing the north pole of $S^n$ under the
stereographical projection, are to consider as k-spheres too.
From the
euclidean point of view the $k$-planes have infinite radius, and
their
centres are not defined; both, centre and radius, are concepts of
euclidean (or Riemannian), but not of M\"obius geometry.
Obviously,
the models differ by the underlying Riemannian metrics only.

\beglem The hypersphere $S(z, r)\subset E^n$ with centre $z$ and
radius $r$ is the 
image under the stereographical projection \em (22) \em of the
hypersphere 
$S(\ffr n)$ defined by \em (18), \em where $\ffr n$ denotes the
spacelike unit
vector
 \begeq \ffr n(z,r) = (2z + \ffr a_{n+1}(-1 -r^2 + |z|^2) + \ffr
a_{n+2}(1 -r^2 + |z|^2))/2r,\:  |z|^2 = \sum_{j=1}^n z_j^2.
 \eneq 
\enlem 

For the proof we first define:

\begdef A set or sequence of $k +2$ points $p_j\in S^n$,
$j=1,\ldots,k+2$, $k\leq n$, is said to be \em in general
position, \em if they do not belong to any $l$-sphere of lower
dimension $l<k$.
\endef
Now we take $n+1$ points of the hypersphere $S(z,r)$ in general
position. Their images under the \em inverse stereographical
projection\em
\begeq
y\in E^n\mapsto sp^{-1}(y) = \frac{2y+\ffr
a_{n+1}(<y,y>-1)}{<y,y>+1}
\eneq 
are $n+1$ points $p_i$ in general position, spanning the
hypersphere
$sp^{-1}(S(z,r))\subset S^n$. The corresponding isotropic vectors
are
linearly independent, and their normed generalised cross product
gives
$\ffr n(z,r)$. The elementary calculation following this way
leads to
formula (24). Finally, calculating the scalar product of two
vectors
of the shape (24) we obtain (21).  \bst

{\sc Remark} 1. Orientation. Changing the order in the sequence
$$(p_1,\ldots,p_{n+1})$$ changes the sign of $\ffr n$ according
to the
signature of the permutation. Therefore $H^{n+1}$ corresponds
bijectively to the manifold of all oriented hyperspheres of the
M\"obius space $S^n$.

{\sc Remark} 2. Hyperplanes. Formulas (21) and (24) are related
to the
euclidean model; they make no sense for hyperplanes. If the
hyperplane
of $E^n$ is defined by the equation $<\ffr b, \ffr x> = p$,
instead of
(24) we get a corresponding vector $\ffr n(\ffr b, p) \in
H^{n+1}$ by
\begeq
\ffr n(\ffr b, p) = \frac{\ffr b + (\ffr a_{n+1} + \ffr
a_{n+2})p}{|\ffr b|},
\eneq 
where $|\ffr b|$ denotes the norm of $\ffr b$. If another
hyperplane
is defined by $\hat{\ffr b}, \hat{p}$, we obtain $$<\ffr n(\ffr
b, p),
\ffr n(\hat{\ffr b}, \hat{p})> = <\ffr b, \hat{\ffr b}>.$$ Here
for
simplicity $\ffr b, \hat{\ffr b}$ are assumed to be unit
vectors. Therefore the conformal invariant is the cosinus of the
angle
between the normal vectors of the intersecting or parallel
hyperplanes; parallel means tangent to each other at infinity. Of
course, this invariant never can be larger than one. Finally, in
case
of the mutual position of a hyperplane and a hypersphere, we have
for
its conformal invariant in euclidean terms
\begeq
<\ffr n(\ffr b, p), \ffr n(z,r)> = \frac{<\ffr b, z> - p}{r|\ffr
b|},
\eneq 
and this coincides with the cosinus of the intersection angle if
the
intersection is not empty.

\section{Pairs of subspheres} 
In this section we classify pairs $(S_1^m, S_2^l)$ of subspheres
of
$S^n$ under the action of the M\"obius group $G_n$. We assume

A. The dimensions fulfil $n>m\geq l\geq 0.$

B. The pair is in \em general position, \em i.e. there does not
exist a 
hypersphere  $S^{n-1}$ with $S_1^m \cup S_2^l\subset S^{n-1}$. 

If assumption B is not fulfilled one considers the smallest
subsphere containing $S_1^m \cup S_2^l$ and applies the methods
described below. We often omit the dimension superscripts and fix
the notations as follows: The pseudo-euclidean subspace of
dimension $m+2$ defining $S_1$ is $\bbm W\in X_m$, and the
$(l+2)$-dimensional pseudo-euclidean subspace defining $S_2$ is
$\bbm U\in X_l$ (see proposition 1). Assumption B is equivalent
to each of the following equations:
\begeq
\bbm W + \bbm U = \bbm V,\; \bbm W^\bot \cap \bbm U^\bot = \ffr
o.
\eneq 
Therefore the sum $\bbm W^\bot +\bbm U^\bot$ is direct. From (15)
and (8) we get
\begin{eqnarray}
&&\dim \bbm U \cap \bbm W = m + l - n + 2 \geq 0,\\
&&\dim \bbm U^\bot + \bbm W^\bot = 2n-m-l. 
\end{eqnarray} 
Relating qualitative different  positions of the subspheres and
algebraic properties of $\bbm U \cap \bbm W$, we shall
distinguish four cases:

Case 1.  $\bbm U \cap \bbm W$ is euclidean, or the null vector
$\ffr o$. Then  $\bbm U^\bot + \bbm W^\bot$ is pseudo-euclidean. 
The subspheres $S_1, S_2$ are disjoint, and there exist
hyperspheres separating them:
\begeq
\Sigma_1 = S(\ffr n)\supset S_1, \Sigma_2 = S(\ffr m) \supset S_2
\mbox{ with }\Sigma_1\cap\Sigma_2 = \empt.
\eneq 
Indeed, by (30) we have $\dim \bbm U^\bot + \bbm W^\bot\geq 2$.
Since $\bbm U^\bot + \bbm W^\bot$  is pseudo-euclidean, one can
find a timelike vector $\ffr x = \ffr n + \ffr m \in \bbm U^\bot
+ \bbm W^\bot$, where $\ffr m \in \bbm U^\bot,\, \ffr n \in \bbm
W^\bot$ must be spacelike. Since their span $[\ffr m, \ffr n]$ is
pseudo-euclidean, the corresponding hyperspheres are disjoint
(section 3, case 2).

Case 2. $\bbm U \cap \bbm W$ is pseudo-euclidean, and we have
$k=\dim \bbm U \cap \bbm W \geq 2 $. Then this intersection
contains isotropic vectors, and defines a subsphere of dimension
$k-2$. 

Case 3. $\bbm U \cap \bbm W$ is pseudo-euclidean, and we have
$k=\dim \bbm U \cap \bbm W = 1 $. Then this one-dimensional,
timelike intersection does not contain isotropic vectors, and the
subspheres do not intersect: $S_1\cap S_2 = \empt$. In difference
to case 1  there do not exist hyperspheres $\Sigma_1, \Sigma_2$
which separate the subspheres, i.e. have properties (31); the
subspheres $S_1, S_2$ are \em interlaced. \em Indeed, for any two
linearly independent vectors $\ffr m, \ffr n$ of the euclidean
subspace $(\bbm U\cap \bbm W)^\bot = \bbm U^\bot + \bbm W^\bot$
the span $[\ffr m, \ffr n]$ is euclidean, and by section 3, case
1, the corresponding hyperspheres intersect. As examples we
mention a hypersphere $S_1$, and a\linebreak 0-sphere $S_2$, the
two points of which lie on different sides of $S_1$, or two
interlaced circles in $S^3$. By (29),  in case 3 we always have
$m+l=n-1$.

Case 4. $\bbm U \cap \bbm W$ is isotropic. Since in any isotropic
subspace of the pseudo-euclidean space $\bbm V$ there exists a
uniquely defined isotropic subspace of dimension one, the
intersection of the subspheres is a uniquely defined point:
$S_1\cap S_2 = \{x_0\}$. Examples are a hypersphere $S_1$ and any
tangent $l$-sphere $S_2$, or two circles in $S^3$ which intersect
in exactly one point; circles in $S^3$ intersecting in two
points, or having a common tangent, are not in general position.

To obtain a complete system of invariants for pairs of subspheres
we apply the method of stationary angles, see e.g.  {\sc
B.~A.~Rosenfeld} \cite{Ro:x}, \S~3.3, \S~11.3,  {\sc
H.~Reichardt} \cite{Rei:68} \S\ 5.3. We consider the function
$<\ffr n, \ffr m>$ under the conditions
\begeq
\ffr n \in \bbm W^\bot, \ffr m \in \bbm U^\bot, <\ffr n, \ffr n>
= 1, <\ffr m, \ffr m> = 1.
\eneq 
Since $\bbm W^\bot, \bbm U^\bot$ are euclidean vector spaces, the
vector pair $(\ffr n, \ffr m)$ varies on a compact set, and
therefore exist maxima and minima of $<\ffr n, \ffr m>$.
Introducing the Lagrange multipliers $\lambda, \mu$ we consider
the function 
\begeq
f(\ffr n, \ffr m, \lambda, \mu) =
<\ffr n, \ffr m> -\lambda(<\ffr n, \ffr n> - 1) - \mu(<\ffr m,
\ffr m> - 1).
\eneq 
At an extremum the differentials of $f$ with respect to $\ffr n,
\ffr m, \lambda, \mu$ must vanish. Deriving  $f$ with respect to
$\ffr n, \ffr m$ we obtain the conditions
\begin{eqnarray}
&&<d\ffr n, \ffr m> - 2\lambda<d\ffr n, \ffr n> = <d\ffr n, \ffr
m- 2\lambda\ffr n> = 0,\\
&&<d\ffr m, \ffr n> - 2\mu<d\ffr m, \ffr m> = <d\ffr m, \ffr n-
2\mu\ffr m> = 0.
\end{eqnarray}
The derivations with respect to $\lambda, \mu$ only reproduce the
conditions (32). The decompositions
\begeq 
\bbm V =\bbm W\oplus\bbm W^\bot= \bbm U\oplus\bbm U^\bot
\eneq 
allow to introduce the linear maps
\begeq
p_{1,2}: \ffr m \in \bbm U^\bot \mapsto pr_{\bbm W^\bot}(\ffr
m)\in \bbm W^\bot,\:
p_{2,1}: \ffr n \in \bbm W^\bot \mapsto pr_{\bbm U^\bot}(\ffr
n)\in \bbm U^\bot.
\eneq 
Keeping in mind the conditions (32) we get from the second
equations (34), (35) the equivalent conditions 
\begeq
 2\lambda \ffr n = p_{1,2}(\ffr m), \: 2\mu \ffr m = p_{2,1}(\ffr
n).
\eneq 
Now we have to find real numbers $\lambda, \mu$ such that this
equations admit non trivial solutions $\ffr n, \ffr m$. We
consider the linear endomorphism
\begeq
A = p_{1,2}\circ p_{2, 1} \in \Eend(\bbm W^\bot).
\eneq 
If  $\ffr n, \ffr m$ is a solution of (38) then it follows 
$$A(\ffr n) = 4\lambda\mu\ffr n.$$
Therefore we shall look for eigenvalues of $A$. The linear maps
$p_{1,2},p_{2,1}, A$ operate between euclidean vector spaces. We
show
\beglem The operator $A\in \Eend(\bbm V)$ is selfadjoint, and the
formulas
\begeq
p'_{2,1} = p_{1,2},\: p'_{1,2} = p_{2,1},\: A' = A =p'_{2,1}\circ
p_{2,1}
\eneq 
are valid. The eigenvalues of $A$ are nonnegative. The eigenspace
of $A$ to the eigenvalue zero is the kernel of $p_{2,1}$.
\enlem 
\pro For arbitrary vectors $\ffr x \in \bbm U^\bot, \ffr y \in
\bbm W^\bot$ we denote the decompositions with respect to the
orthogonal direct sums (36) by
$$
\ffr x = \ffr x_{W} + \ffr x_{W^\bot}, 
\ffr y = \ffr y_{U} + \ffr y_{U^\bot}.
$$
By the definitions we have
$$<p'_{2,1}(\ffr x),\ffr y> = <\ffr x, p_{2,1}(\ffr y)> = <\ffr
x, \ffr y_{U^\bot}> = <\ffr x, \ffr y>.$$
The last equation follows because $\ffr x \in \bbm U^\bot$ is
orthogonal to $\ffr y_{U}$. Since $\ffr y\in\bbm W^\bot$ we
continue
$$ <\ffr x, \ffr y>= <\ffr x_{W^\bot}, \ffr y> = <p_{1,2}(\ffr
x), \ffr y>.$$
Comparing the first and the last term of these equations we
obtain the first statement of (40), from which the other both
follow easily. Now let $\ffr a$ be an eigenvector of $A$ to the
eigenvalue $\alpha$. By (40) it follows
$$<\ffr a, A(\ffr a)> = \alpha <\ffr a, \ffr a> = <p_{2,1}(\ffr
a),p_{2,1}(\ffr a)>.$$
Since $\bbm W^\bot, \bbm U^\bot$ are euclidean vector spaces, and
$\ffr a\neq\ffr o$ we obtain $\alpha\geq 0$, and the last
statement. \bst

Because $A$ is a  selfadjoint endomorphism, it has $n-m$ real
eigenvalues
\begeq
\alpha_1\geq \alpha_2 \geq \ldots \geq \alpha_{n-m}\geq 0.
\eneq 
Let $(\ffr a_i),\, i = 1,\ldots , n-m,$ be an orthonormal basis
of eigenvectors with $A(\ffr a_i) =\alpha_i \ffr a_i$. Applying
(40) we obtain
\begeq
<p_{2,1}(\ffr a_i),p_{2,1}(\ffr a_j)> = <\ffr a_i, A(\ffr a_j)> =
<\ffr a_i, \ffr a_j> \alpha_j = \delta_{i\, j}\alpha_j.
\eneq 
Let $r_A$ denote the rank of $A$. We define
\begeq
\ffr b_j := p_{2,1}(\ffr a_j)/\sqrt{\alpha_j},\: j = 1,\ldots,
r_A.
\eneq 
 Applying (40) one easily calculates
\begin{eqnarray}
&& <\ffr b_i, \ffr b_j> = \delta_{i\, j},\: i, j = 1,\ldots,r_A\\
&& <\ffr b_i, \ffr a_k> = \delta_{i\, k}\sqrt{\alpha_k},\: k =
1,\ldots, n-m.
\end{eqnarray}
The extremal properties of the eigenvalues show that the angles
defined by 
\begeq       
\cos \beta_k = \sqrt{\alpha_k} \mbox{, if } \alpha_k \leq 1,0\leq
\beta_k\leq \pi/2,
\eneq 
can be considered as the stationary angles between the
hyperspheres fulfilling (31). Of course, this is true only if the
hyperspheres intersect, as we discussed in section 3. Generally
one can speak only about the stationary values of $<\ffr n, \ffr
m>$ under the conditions (31). Now we complete the orthonormal
sequence $(\ffr b_i), i = 1,\ldots, r_A,$ to an orthonormal basis
of $\bbm U^\bot$. As a consequence of Lemma 4 we have
$$ <p_{1,2}(\ffr b_i), \ffr a_j>=<\ffr b_i, p_{2,1}(\ffr a_j)> =
<\ffr b_i,\ffr b_j>\sqrt{\alpha_j} = \delta_{i\,
j}\sqrt{\alpha_j}$$
for $i = 1,\ldots, n-l, j = 1, \ldots, n-m$. If
 $j\leq r_A$ this follows from  definition (43), and for $j>r_A$
we remember that by lemma 4 we have
\begeq
\ker p_{2,1} = \bbm W^\bot \cap \bbm U = [\ffr a_{r_A+1},\ldots,
\ffr a_{n-m}].
\eneq 
It follows
\begeq
p_{1,2}(\ffr b_j) = \ffr a_j\sqrt{\alpha_j}\mbox{ for } j\leq 
r_A,\; p_{1,2}(\ffr b_k) = \ffr o \mbox{ for } k>r_A.
\eneq 
Thus, in analogy with (47) we conclude
\begeq
\ker p_{1,2} = \bbm U^\bot \cap \bbm W = [\ffr b_{r_A+1},\ldots,
\ffr b_{n-l}].
\eneq 

For the following the maximum $\alpha_1$ of the eigenvalues is
deciding. We shall relate its value to the four cases discussed
at the beginning of this section. We define 
\begdef
 The orthogonal bases $(\ffr a_i),\, (\ffr b_j)$ of $\bbm V$
fulfilling (1)  are said to be \em adapted to the subspheres
$S_1, S_2$,\em\  if \\ 
1. $(\ffr a_1,\ldots, \ffr a_{n-m})$ is a  basis of eigenvectors
of $A$ such that (41) and $A(\ffr a_i) = \ffr a_i\alpha_i$ are
valid;\\
2. $(\ffr b_1,\ldots, \ffr b_{n-l})$ is a basis of $\bbm U^\bot$
such that (45) is satisfied for $i = 1,\ldots,n-l$. \endef

\beglem
Under the assumptions \em  A., B. \em the following disjunction
is valid:
\begin{eqnarray*}
\alpha_1 > 1 & \Longleftrightarrow & \bbm U \cap \bbm W \mbox{
euclidean, }\\
\alpha_1 = 1 & \Longleftrightarrow & \bbm U \cap \bbm W \mbox{
isotropic, } \\
\alpha_1 < 1 & \Longleftrightarrow & \bbm U \cap \bbm W \mbox{
pseudo-euclidean }
\end{eqnarray*}
\enlem 
\pro We consider $\bbm U^\bot + \bbm W^\bot = (\bbm U\cap\bbm
W)^\bot$, and adapt the bases. If $\alpha_1 > 1$, we conclude
from (45) that the span $[\ffr a_1, \ffr b_1]$ is a
pseudo-euclidean subspace of  $\bbm U^\bot + \bbm W^\bot$, which
therefore must be pseudo-euclidean too. Conversely, if this is
the case, we find a timelike vector $\ffr x = \ffr u + \ffr w$,
$\ffr u\in\bbm U^\bot, \ffr w\in\bbm W^\bot$. Let $\ffr u_0, \ffr
w_0$ be the corresponding normed vectors. Since the span $[\ffr
u, \ffr w]=[\ffr u_0, \ffr w_0]$ is pseudo-euclidean, the
determinant of the scalar products $1 - <\ffr u_0, \ffr w_0>^2$
must be negative, and by the maximum property of the eigenvalue
$\alpha_1$ we get
$$1 <\; <\ffr u_0, \ffr w_0>^2 \:\leq \alpha_1.$$
Now let $\alpha_1 = 1$. By (45) we have $<\ffr b_1, \ffr a_1> =
1$. Therefore the vector $\ffr z =\ffr a_1 -\ffr b_1$ satisfies
$<\ffr z, \ffr z> = 0$. By (28) we have $\ffr z\neq \ffr o$, and
therefore $\ffr z$ is isotropic. From  (45) it  follows $<\ffr z,
\ffr a_k> = 0$ for $k = 1,\ldots, n-m$, and $<\ffr z, \ffr b_i> =
0$ for $i = 1,\ldots, r_A$. Finally, from (49) we obtain  $<\ffr
z, \ffr b_j> = 0$ for $j = r_A + 1,\ldots, n - l$. Therefore 
$\bbm U^\bot + \bbm W^\bot$ is isotropic, and the orthogonal
space  $\bbm U\cap \bbm W$ too. In the converse, assume that 
$\bbm U^\bot + \bbm W^\bot$ is isotropic. Then we find a vector
$\ffr y = \ffr u + \ffr w\neq \ffr o$ being orthogonal to  $\bbm
U^\bot + \bbm W^\bot$, with $\ffr u\in\bbm U^\bot, \ffr w\in\bbm
W^\bot$. Especially we have
$$<\ffr y, \ffr u> = <\ffr u, \ffr u> + <\ffr w, \ffr u> = 0,\;
<\ffr y, \ffr w> = <\ffr u,\ffr w> + <\ffr w, \ffr w> = 0,$$
which yields 
$$\frac{ <\ffr w, \ffr u>^2}{ <\ffr u, \ffr u> <\ffr w, \ffr w>}
= 1.$$
By the maximum property of $\alpha_1$ we conclude $\alpha_1 \geq
1$, but $\alpha_1 > 1$ cannot take place since then $\bbm U^\bot
+ \bbm W^\bot$ would be pseudo-euclidean, as already shown. The
third equivalence is an immediate consequence of the others and
(41). \bst

By (43) and (47) we know the components of the eigenvectors of
$A$ in $\bbm U^\bot$. Our aim is now to calculate the components
of these vectors in $\bbm U$. Let $q: \bbm W^\bot \to \bbm U$
denote the projection defined by the orthogonal decompositions
(36). Using $q = \id_{\bbm W^\bot} - p_{2,1}$ and applying (43) -
(45) one calculates
\begeq
<q(\ffr a_i),q(\ffr a_j)> = \delta_{i\, j}(1 -  \alpha_j).
\eneq 

Case 1. $\alpha_1 > 1$. \em Then the space $\bbm U^\bot + \bbm
W^\bot$ is pseudo-euclidean,  the spheres are disjoint: $S_1\cap
S_2 =\empt$, and can be separated by hyperspheres. The eigenvalue
$\alpha_1$ has multiplicity \em 1\em , and the eigenvalues
$\alpha_i, i\geq 2,$ satisfy $\alpha_i < 1$. We remember \em
(46)\em , and set 
$$ \sqrt{\alpha_1} = \cosh \beta_1,\; \beta_1 > 0.$$
There exist adapted bases such that \em
\begin{eqnarray}
&&\ffr a_1 = \ffr b_1 \cosh \beta_1 + \ffr b_{n+2}\sinh \beta_1\\
&&\ffr a_i = \ffr b_i \cos \beta_i + \ffr b_{n-l+i -1} \sin
\beta_i,\; (i=2,\ldots,n-m).
\end{eqnarray} 
\pro By Lemma 5 the space $\bbm U^\bot + \bbm W^\bot$ is 
pseudo-euclidean. By norming $q(\ffr a_1)$ we define $\ffr
b_{n+2}$  and obtain (51).  There can not be any other
eigenvector with an eigenvalue $\alpha > 1$, since then, by (50), 
we would have two orthogonal timelike vectors, what is impossible
in a pseudo-euclidean space of index 1. In particular, $\alpha_1$
has multiplicity 1. Assume $\alpha_2 = 1$. Then from (50), for
$i=j=2$, would follow $q(\ffr a_2) = \ffr o$, or $q(\ffr a_2)$
isotropic. The latter cannot take place, since this vector is
orthogonal to $\ffr b_{n+2}$ and therefore belongs to an
euclidean subspace. On the other hand, $q(\ffr a_2) = \ffr o$
implies $\ffr a_2 \in \bbm W^\bot\cap\bbm U^\bot$, what
contradicts (28). Thus, $0\leq \alpha_i < 1$ holds true for $i =
2,\ldots, n-m$. Norming the vectors $q(\ffr a_i)$ we find the
orthonormal sequence $(\ffr b_{n-l+i-1})_i$ such that (52) holds.
Completing the already defined vectors to orthonormal bases of
$\bbm U$ and $\bbm U^\bot$, respectively, we get the adapted
bases with properties (51), (52).\bst

Case 2. \em  $\alpha_1 < 1$, and  $\dim \bbm U\cap \bbm W > 1$.
Then, by \em (29)\em\  and Lemma 5, $\bbm U\cap\bbm W$ is
pseudo-euclidean, and the intersection $S_1\cap S_2$ is a
subsphere of dimension $m+l-n$. There exist adapted bases such
that\em 
\begeq
\ffr a_i = \ffr b_i \cos \beta_i + \ffr b_{n-l+i} \sin \beta_i,\;
(i=1,\ldots,n-m).
\eneq 
\pro  Equation (46) is a correct definition for $k = 1, \ldots,
n-m$, giving stationary angles of the intersecting subspheres. 
By (50) we obtain as projections in $\bbm U$ an orthogonal
sequence, which we can norm and take as part of the adapted
basis. Analogously to (52) we obtain (53). \bst

Case 3. \em  $\alpha_1 < 1$, and $\dim \bbm U\cap\bbm W = 1$.
Then  the subspheres are interlaced and one has $m+l-n = -1$.
There exist  adapted bases  such that \em (53) \em\ holds.\em\\
\pro The first statement is already proved, see case 3 at the
beginning of this section and (29).The rest can be proved  as in
case 2. \bst

Case 4. \em  $\alpha_1 = 1$. Then the multiplicity of this
eigenvalue is 1. The intersection $\bbm U\cap \bbm W$ is
isotropic, and the subspheres intersect in a single point.  There
exist $n-m-1$ stationary angles defined by \em (46)\em\  for $k =
2,\ldots, n-m$ and adapted bases such that \em (52) \em and the
following equation hold:\em
\begeq
\ffr a_1 = \ffr b_1 + \ffr b_{n+1}+\ffr b_{n+2}.
\eneq 
\pro  By Lemma 5 the second statement is true. Since $q(\ffr
a_1)$ is an  isotropic vector in $\bbm U$ which by (50) is
orthogonal to $q(\bbm W^\bot)$ we find an orthogonal  basis
$(\ffr b_j)$, $j = n-l+1,\ldots,n+2$ of $\bbm U$, such that (54)
holds. If the multiplicity of $\alpha_1 = 1$ would be greater
than 1, then an analogous representation would exist for $\ffr
a_2$, where the isotropic parts were proportional, since the
one-dimensional isotropic subspace in $q(\bbm W^\bot)$ is
uniquely defined. Thus, together with (54) we would have
$$\ffr a_2 = \ffr b_2 + (\ffr b_{n+1}+\ffr b_{n+2})\kappa$$
for a certain constant $\kappa\neq 0$. Multiplying (54) by
$\kappa$, and subtracting, we would get a non-vanishing vector
$$\ffr a_1\kappa -\ffr a_2 =\ffr b_1\kappa - \ffr b_2 \in \bbm
U^\bot\cap\bbm W^\bot,$$
in contradiction to (28). Now, for $k=2,\ldots, n-m$, we may
apply (46); norming and numbering the vectors $q(\ffr a_k)$ in an
appropriate way, we obtain (52). Since $q(\ffr a_1)$ is
orthogonal to $q(\bbm W^\bot)$, the orthogonal complement of the
span $[q(\ffr a_2),\ldots,q(\ffr a_{n-m}]$ in $\bbm U$ is an at
least two-dimensional pseudo-euclidean subspace containing
$q(\ffr a_1)$, in which we may realize the equation (54);
completing to orthogonal bases of $\bbm U, \bbm U^\bot$ we get
the required adapted bases. \bst

We summarise the results and finish the article with the
following

\begin{theo}
Under the assumptions A, B the eigenvalues \em (41)\em\  of the
selfadjoint operator $A$, defined by \em (39),\em\  are a
complete system of invariants under the M\"obius group for pairs
of subspheres $(S_1^m, S_2^l)$ of the M\"obius space $S^n$. The
cases \em\ 1 - 4\em\  of the mutual position of the spheres
correspond to the cases\em\  1 - 4\em\   characterised by the
maximal eigenvalue.
\end{theo}
\pro Obviously, by their definition, the eigenvalues $\alpha_i$
are invariant under M\"obius transformations. In each of the
cases we proved the existence of adapted bases with respect to
which the eigenbasis $(\ffr a_i)$ of the space $\bbm W^\bot$ is
expressed with respected to the basis $(\ffr b_i)$, $i =
1,\ldots, n+2$, corresponding to the subsphere $S^l_2$, with
coefficients uniquely defined by the eigenvalues $\alpha_i$.
Therefore, if another pair $(\hat{S}_1^m, \hat{S}_2^l)$ possesses
the same eigenvalues, and $(\hat{\ffr a}_i),\, (\hat{\ffr b}_i)$
are the corresponding  adapted bases, then the M\"obius transform
$g$ defined by $g(\ffr b_i) = (\hat{\ffr b}_i)$, transforms the
pairs into each other:
$$(gS_1,\, gS_2) = (\hat{S}_1,\, \hat{S}_2).$$
We remark that case 2 and case 3 differ by the dimensions: In
case 3 we have $ m+l-n = -1$, and in case 2 $m + l-n\geq 0$ is
the dimension of the subsphere $S_1\cap S_2$. \bst
\pagebreak

\bibliographystyle{plain}

\vspace{15mm}

 e-mail: Rolf.Sulanke@t-online.de 

\end{document}